\documentclass[letterpaper, 12 pt]{article}
\usepackage{amsthm, amsmath, amssymb, amsfonts}

\newcommand{\setleftmargin}[1]{
        \addtolength{\textwidth}{\oddsidemargin}
        \addtolength{\textwidth}{1in}
        \addtolength{\textwidth}{-#1}
        \setlength{\oddsidemargin}{-1in}
        \addtolength{\oddsidemargin}{#1}
        \setlength{\evensidemargin}{\oddsidemargin}
}
\newcommand{\setrightmargin}[1]{
        \setlength{\textwidth}{8.5in}
        \addtolength{\textwidth}{-\oddsidemargin}
        \addtolength{\textwidth}{-1in}
        \addtolength{\textwidth}{-#1}
}
\setleftmargin{1 in}
\setrightmargin{1 in}

\newcommand{\Edge}{\mathrm{Edge}}
\newcommand{\Vertices}{\mathrm{Vert}}
\newcommand{\MARKED}{\textbf{MARKED}}
\newcommand{\TRUE}{\textbf{TRUE}}
\newcommand{\FALSE}{\textbf{FALSE}}
\newcommand{\into}{\hookrightarrow}

\newcommand{\sym}{\sim}

\newtheorem{conj}{Conjecture}
\newtheorem{Theorem}[conj]{Theorem}
\newtheorem{prop}[conj]{Proposition}

\newtheorem{Lemma}[conj]{Lemma}

\title{Cyclically Orientable Graphs}
\author{David E Speyer                }

\begin{document}

\maketitle

\abstract{Barot, Geiss and Zelevinsky define a notion of a
  ``cyclically orientable graph'' and use it to devise a test for
  whether a cluster algebra is of finite type. Barot, Geiss and
  Zelivinsky's work leaves open the question of giving an efficient
  characterization of cyclically orientable graphs. In this paper, we
  give a simple recursive description of cyclically orientable graphs,
  and use this to give an $O(n)$ algorithm to test whether a graph on
  $n$ vertices is cyclically orientable. Shortly after writing this
  paper, I learned that most of its results had been obtained
  independently by Gurvich \cite{Gurv}; I am placing this paper on the arXiv
  to spread knowledge of these results.}

\section{Introduction and Results}

In ``Cluster Algebras and Positive Matrices'' \cite{BGZ}, a graph $G$
is defined to be cyclically orientable if it has an orientation in which
every cycle of $G$ which occurs as an induced subgraph is cyclically
oriented.  The aim of this note is to prove the following characterization of
cyclically orientable graphs:

\begin{Theorem} \label{Main}
  A graph $G$ is cyclically orientable if and only if all of its
  two-connected components are. A two-connected graph is cyclically
  orienteable if and only if it is either a cycle, a single edge, or
  of the form $G' \cup C$ where $G'$ is a cyclically orientable graph,
  $C$ is a cycle and $G'$ and $C$ meet along a single edge. Moreover,
  if $G=G' \cup C$ is any such decomposition of $G$ into a cycle and a
  subgraph meeting along a single edge, then $G$ is cyclically
  orientable if and only if $G'$ is.
\end{Theorem}

It follows easily from this characterization that every two-connected
cyclically orientable graph is series-parallel. In particular,
cyclically orientable graphs are always planar and have at most $O(n)$
edges, where $n$ is the number of vertices.

We use our results to give an $O(n)$ algorithm to determine whether a
graph on $n$ vertices is cyclically orientable or not. Using this
algorithm, it is easy to efficiently test condition (4) of \cite{BGZ}:
the algorithm in this paper not only tests whether or not a graph is
cyclically orientable but also, if that graph is orientable, will find
such an orientation. By Propositions 1.4 and 1.5 of that paper, once
we know how to find a cyclic orientation of a given graph, testing
condition (4) amounts to simply checking whether a certain symmetrizable
matrix is positive definite, which is a standard linear algebra computation.

Shortly after writing this paper, I learned Gurvich had independently
found its main result approximately nine months earlier. With his
consent, I am placing this paper online so that this result will
become known and available to those wising to compute with cluster
algebras. 

\section{Graph Theoretic Terminology}

A graph $G$ is a finite set $\Vertices(G)$ and a subset $\Edge(G)$ of the
set of two element subsets of $\Vertices(G)$. $\Vertices(G)$ and $\Edge(G)$
are called the \emph{vertices} and \emph{edges} of $G$; the elements
of an edge of $G$ are called the \emph{endpoints} of that edge. 
A subgraph of $G$ is a graph $G'$ equipped with injections $\Vertices(G')
\into \Vertices(G)$, $\Edge(G') \into \Edge(G)$ compatible with
containment. If $S$ is a subset of $\Vertices(G)$, $G|_S$ is the subgraph
whose vertices are $S$ and whose edges are the elements of $\Edge(G)$
that are subsets of $S$. A graph of the form $G|_S$ is called \emph{an
  induced subgraph of $D$}. We write $G \setminus S$ for $G|_{\Vertices(G) \setminus S}$.

A \emph{path} in $G$ is a subgraph isomorphic to the graph on the
vertex set $\{ 1,\ldots, n \}$ whose edges are $\{ 1,2 \}$, \ldots, $
\{ n-1, n \}$. The vertices corresponding to $1$ and $n$ are called
the \emph{endpoints} of the path.  A \emph{cycle} in $G$ is a subgraph
isomorphic to the graph on the vertex set $\{ 1,\ldots, n \}$ whose
edges are $\{ 1,2 \}$, \ldots, $ \{ n-1, n \}$, $\{ n,1 \}$. A path is
called a \emph{chain} if it occurs as an induced subgraph; a cycle is
called a \emph{chordless cycle} if it occurs as an induced subgraph.

An orientation of a graph $G$ is an assignment of an ordering of the
endpoints of each edge of $G$. An orientation of a cycle is called
cyclic if it recieves the orientation $(1,2)$, \ldots, $(n-1,n)$,
$(n,1)$ or the opposite orientation. An orientation is called cyclic
if its restriction to every chordless cycle is cyclic. A graph is
called cyclically orientable if it has a cyclic orientation. Note that,
if $G$ is cyclically orientable, so is $G|_S$ for any $S \subseteq
\Vertices(G)$.

Define an equivalence relation on the vertices of $G$ by setting $v_1
\sym v_2$ if there is a path in $G$ with endpoints $v_1$ and $v_2$.
The subgraphs of $G$ induced by the equivalence classes of $G$ are
caled the \emph{connected components} of $G$. $G$ is called
\emph{connected} if it has only one connected component.

Define an equivalence relation on the edges of $G$ by setting $e_1
\sym_2 e_2$ if there is a cycle containing $e_1$ and $e_2$. It is not
obvious, but it is true, that this defines an equivalence relation. It
turns out that the equivalence classes for $\sym_2$ are the edge sets
of unique connected induced subgraphs of $G$. Define the
\emph{two-connected components} of $G$ to be these subgraphs. See chapter III of \cite{Tutte} for background on
two-connectivity.

\section{A Decomposition Result}

The aim of this section is to prove the following result:

\begin{Theorem} \label{Key}
  Let $G$ be a two-connected, cyclically oriented graph which is not a
  cycle or a single edge. Then there exists an edge $e=\{v,w\}$ of $G$
  such that $G \setminus \{ v,w \}$ is disconnected.
\end{Theorem}

\begin{Lemma} \label{Chain}
  Let $G$ be a graph, $v$ and $w$ vertices of $G$ and $P \subseteq
  \Vertices(G)$ the vertices of a path with endpoints $v$ and $w$. Then
  there is a chain joining $v$ and $w$ whose vertices are contained in
  the vertices of $P$.
\end{Lemma}

\begin{proof}
  Let $C$ be the shortest path from $v$ to $w$ whose vertices are
  contined in those of $P$. Let $v=c_1$, $c_2$, \dots, $c_{\ell-1}$,
  $c_{\ell}=w$ be the vertices of $C$. If there were an edge between
  $c_i$ and $c_j$ for $j>i+1$, then $c_1$, $c_2$, \dots, $c_i$, $c_j$,
  \dots, $c_{\ell}$ would be a shorter path, a contradiction. Thus, no
  such edge exists and $C$ is a chain.
\end{proof}

\begin{Lemma} \label{Chordless}
  Let $G$ be a graph, $v$ a vertex of $G$ and $Z \subseteq
  \Vertices(G)$ the vertices of a cycle containing $v$. Then
  there is a chordless cycle containing $v$ whose vertices are contained in
  the vertices of $Z$.
\end{Lemma}

\begin{proof}
Similar to the previous lemma.
\end{proof}

\begin{proof}[Proof of theorem \ref{Key}]
  Fix a cyclic orientation of $G$. Since $G$ is two-connected and not
  a single edge, it does not contain any vertices of degree $1$.
  Moreover, since it is connected and not a single cycle, not all of
  its vertices can be of degree $2$. Thus, $G$ has a vertex of degree
  greater than or equal to $3$, call it $v$. Let $N \subseteq
  \Vertices(G)$ be the neighbors of $v$.
  
  We define a graph $\Gamma$ whose vertices are $N$ and for which
  there is an edge between $u$ and $u' \in N$ if and only if there is
  a path from $u$ to $u'$ in $G \setminus \{ v \}$ whose internal
  vertices are not in $N$. We claim that $\Gamma$ is connected. Proof:
  two vertices lie in the same connected component of $\Gamma$ if and
  only if they lie in the same connected component of $G \setminus \{
  v \}$. By assumption, $G$ is two-connected, so $G \setminus \{ v \}$
  is connected.
  
  We color the vertices of $\Gamma$ black and white; $u$ is
  colored black if $\{ u,v \}$ is oriented towards $v$ and white if it
  is oriented away from $v$. We claim that $\Gamma$, with this
  coloring, is bipartite. Proof: assume for contradiction that $u$ and
  $u'$ are two vertices of $\Gamma$ of the same color and $P \subseteq
  G \setminus \{ v \}$ a path between them not passing through any
  other vertex of $N$. Then $v$ and $P$ form a cycle. By Lemma
  \ref{Chordless}, there is a chordless cycle $C$ containing $v$ and
  all of whose other vertices lie in $P$. As $u$ and $u'$ are the only
  vertices of $P$ bordering $v$, they must be the neighbors of $v$ in
  $C$. But then $C$ can not be cyclically oriented, a contradiction.
  
  So, $\Gamma$ is a bipartite connected graph with at least three
  vertices. We will now show that $\Gamma$ is a tree. Since
  $\Gamma$ is connected, if $\Gamma$ is not a tree, then it contains a
  cycle which, by lemma \ref{Chordless} we may take to be chordless;
  let $u_1$, $u_2$, \ldots, $u_{2 k}$ be the vertices of this cycle,
  with $u_{2i}$ white and $u_{2i+1}$ black. Let $P_i$ be a path
  joining $u_{2i}$ to $u_{2i+1}$ in $G \setminus \{ v \}$ and whose interior vertices do not lie in $N$. By lemma \ref{Chain}, we
  may assume that each $P_i$ is a chain. Then, for every $i$, $v$ and
  $P_i$ form a chordless cycle in $G$ and the assumption that $G$ is
  cyclically oriented implies that the chains $P_i$ are all oriented
  towards their black ends. 
  
  We see that the cycle $\bigcup_i P_i$ in $G$ is not cyclically
  oriented, so it must have a chord; let $x$ and $x'$ be the endpoints
  of this chord. $\{ x,x' \}$ can not be of the form $\{ u_i, u_{i+1}
  \}$, as otherwise $P_i$ would be the edge $\{ u_i, u_{i+1} \}$ (the
  $P_i$ are chains) and $\{ x,x' \}$ would not be a chord. $\{ x,x'
  \}$ can not be of the form $\{ u_i, u_j \}$ for $i$ and $j$ not
  consecutive as then $\{ x,x' \}$ would be a chord of the cycle
  $u_1$, $u_2$, \dots $u_{2k}$ in $\Gamma$. So at least one of $x$ and
  $x'$ is not in $N$, say $x'$. Let $x'$ lie in the interior of the
  path $P_{i'}$. $x$ can not also lie in the path $P_{i'}$, as $P_{i'}$ is a chain.
  
  We consider two cases. If $x=u_i$, then $i \neq i'$ and $i \neq
  i'+1$. Let $i'+\epsilon$ be the same color as $i$, where $\epsilon
  \in \{ 0,1 \}$. Then $u_i=x$, $x'$ and the portion of $P_i$ running
  from $x'$ to $u_{i'+\epsilon}$ form a path in $G \setminus \{ v \}$
  in which none of the internal vertices lie in $N$. So $u_i$ and
  $u_{i'+\epsilon}$ are joined by an edge in $\Gamma$, contradicting
  that $\Gamma$ is bipartite.
  
  If $x$ is not one of the $u_i$, let $x \in P_i$; $i \neq i'$ as
  otherwise $P_i$ would not be a chain. Then we can find $\epsilon$
  and $\epsilon' \in \{ 0,1 \}$ such that $i+\epsilon \neq
  i+\epsilon'$ and $i+\epsilon$ and $i'+\epsilon'$ have the same
  color. Then the path from $u_{i+\epsilon}$ to $x$ in $P_i$, the edge
  from $x$ to $x'$ and the path from $x'$ to $u_{i'+\epsilon'}$ in
  $P_{i'}$ is a path from $u_{i+\epsilon}$ to $u_{i'+\epsilon'}$
  contained in $G \setminus v$ and with all of its internal vertices
  not in $N$. So $u_{i+\epsilon}$ and $u_{i'+\epsilon'}$ are joined in
  $\Gamma$, contradicting that $\Gamma$ is bipartite. We have now
  shown that $\Gamma$ is a tree.
  
  $\Gamma$ is a tree with at least three vertices. Therefore, there is
  a vertex $w$ of $\Gamma$ so that $\Gamma \setminus \{ w \}$ is
  disconnected. It is then easy to see that $G \setminus \{ v,w \}$ is
  disconnected.
\end{proof}

\section{Proof of Theorem \ref{Main}}

If $G_1$ and $G_2$ are graphs and $e_1$ and $e_2$ are edges of $G_1$ and $G_2$ respectively,
each equipped with an orientation, then we write $G_1 \cup_{e_1,e_2} G_2$ for the
graph formed by gluing $G_1$ and $G_2$ together along $e_1$ and $e_2$ in a
manner compatible with the orientations of $e_1$ and $e_2$. As a preliminary to proving Theorem \ref{Main}, we show.

\begin{Lemma} \label{Gluing}
  Let $G$ be a two-connected, cyclically orientable graph which is not a cycle or a single edge. Then we can
  write $G$ as $G_1 \cup_{e_1, e_2} G_2$ where each $G_i$ is a
  two-connected, cyclically orientable graph which is not a single edge and $e_i$ is an oriented edge of $G_i$. Moreover, every graph of this form is cyclically orientable.
\end{Lemma}

\begin{proof} 
  Let $e= \{ v,w \}$ be the edge of $G$ found in theorem \ref{Key}, so
  that $G \setminus \{ v,w \}$ is disconnected. Let $\Vertices(G \setminus
  \{ v,w \}) = S_1 \sqcup S_2$, where there are no edges connecting
  $S_1$ to $S_2$ and neither of the $S_i$ is empty. We will take $G_i=G|_{S_i \cup \{ v,w 
\}}$ and $e_i$
  to be the edge $\{ v,w \}$ of $G_i$, oriented from $v$ to $w$.
  Clearly, $G = G_1 \cup_{e_1,e_2} G_2$. Also, as each $G_i$ is an
  induced subgraph of the cyclically orientable graph $G$, the $G_i$ are
  cyclically orientable. As both $S_i$ are nonempty, $G_i$ is not a single edge. Thus, it 
remains to show that the $G_i$ are two-connected.
  
  Let $e'=\{x,y\}$ be an edge of $G_i$. We will show that there is a
  cycle, contained entirely in $G_i$, and containing edges $e$ and
  $e'$. This will show that $G_i$ is two-connected. 
  
  Since $G$ is two-connected, we can find a cycle $C$ in $G$
  containing $e$ and $e'$. We claim that actually this cycle must lie
  entirely in $G_i$. Let $p$, $v$, $w$, $q$ be the four vertices of
  $C$ nearest to $e$. Then $C \setminus \{ v,w \}$ is a path from $p$
    to $q$ contained entirely within $G \setminus \{ v,w \}$ and hence
    either entirely within $G|_{S_1}$ or $G|_{S_2}$. As this path
    contains $e' \in \Edge(S_i)$, it lies entirely in $G|_{S_i}$. Then
    $C$ lies entirely in $G|_{S_i \cup \{ v,w \}}= G_i$, as desired.
    
    For the converse direction, take cyclic orientations of $G_1$ and
    $G_2$; after possibly reversing them we can impose that they agree
    with the given orientations on $e_1$ and $e_2$. We obtain an
    orientation of $G:=G_1 \bigcup_{e_1,e_2} G_2$ in the obvious way.
    Let $e$ denote the edge of $G$ coming from the $e_i$ and $\{ v,w
    \}$ its endpoints.  Let $C$ be any chordless cycle contained in
    $G$. As explained in the previous paragraph, if $C$ contains $e$
    as an edge then $C$ must lie entirely in either $G_1$ or $G_2$ and
    thus must be cyclically oriented. If $C$ does not contain both $v$
    and $w$ then $C$ must again lie solely in $G_1$ or in $G_2$ and be
    cyclically oriented. If $C$ does contain both $v$ and $w$ but not
    consecutively, then $e$ is a chord of $C$, a contradiction. So every chordless cycle of 
$G$ is cyclically oriented. 

\end{proof}

We now prove theorem \ref{Main}.

\begin{proof}
  If all of the two-connected components of a graph are cyclically
  oriented, this clearly provides a cyclic orientation on the graph as
  a whole, because any cycle lies in a single two-connected component.
  A cycle or a single edge is clearly cyclically orientable. If $G'$
  is cyclically orientable and $G=G' \bigcup_{e_1,e_2} C$ for $C$ a
  cycle and $e_1$, $e_2$ oriented edges of $C$ and $G'$ repectfully,
  it follows from the last sentence of lemma \ref{Gluing} that $G$ is
  cyclically orientable. Conversely, if $G=G' \bigcup_{e_1,e_2} C$ and
  $G$ is cyclically orientable then $G'$ is an induced subgraph and
  hence cyclically orientable.  We have now checked all the claims of
  the theorem except for the claim that any two-connected cyclically
  orientable $G$ is of the form $G' \cup_{e_1,e_2} C$.

  Let $G$ be a two-connected cyclically orientable graph. If $G$ is a
  cycle, we are done. If not, use lemma \ref{Gluing} to write $G=G_1
  \bigcup_{e_1,e_2} G_2$, with notation as in that lemma. If $G_2$ is
  a cycle, we are done. If not, by induction on the number of vertices
  in $G$ (and using that $G_2$ is not a single edge), we can write
  $G_2=H \cup_{f_1, f_2} C$ for $C$ a cycle, $H$ a two-connected
  cyclically orientable graph and $f_1$, $f_2$ edges of $H$ and $C$.
  Then
  $$G=G_1 \bigcup_{e_1,e_2} \left( H \bigcup_{f_1,f_2} C \right) =
  \left( G_1 \bigcup_{e_1,e_2} H \right) \bigcup_{f_1,f_2} C.$$
  By the
  last sentence of lemma \ref{Gluing}, $\left( G_1 \bigcup_{e_1,e_2} H
  \right)$ is a two-connected, cyclically orientable graph, so we are
  done.
\end{proof}

\section{The Algorithm}

The purpose of this section is to present an algorithm to test whether
a graph $G$ is cyclically orientable. We first presnt a naive, but
easy to follow, implementation which runs in $O(n^2)$ time, where $n$
is the number of vertices of $G$. We then give a more careful
implementation which runs in $O(n)$ time. For the first result, pretty
much any computation model and presentation of a graph is equivalent
to any other. For our second result we need to assume that we are
using a pointer machine -- that is, a machine which can follow a
pointer into an arbitrarily large memory in a single step -- and that
our graph is presented as a list of vertices with, for each vertex, a
list of pointers from that vertex to its neighbors. It
does not seem likely that there will be a need to do any very large
computations of this sort so precise error bounds are not that
important; we include them because it seems a shame not to point them
out.

We can find the two-connected components of $G$ in $O(e)=O(n^2)$ steps
(see \cite{TwoConn}), so in our first algorithm we can reduce to the
case where $G$ is two-connected. Here is a presentation of our naive
algorithm to test whether a two-connected graph $G$ is cyclically
orientable. Our algorithm also uses a boolean function $\MARKED$ which
assigns the value $\TRUE$ or $\FALSE$ to every degree two vertex of
$G$.

\begin{enumerate}

\item Test whether $G$ is a single edge. If so, return ``YES''.

\item Set $\MARKED(v)=\FALSE$ for every degree two vertex of $G$.
  
\item Find a degree $2$ vertex $v$ of $G$ for which
  $\MARKED(v)=\FALSE$ or determine that none exists ($O(n)$ steps). If
  none exists, return ``NO''.
  
\item Find the unique path $u_1$, $u_2$, \dots, $u_{k-1}$, $u_k=v$,
  $u_{k+1}$, \dots, $u_{\ell}$ in $G$ such that $u_i$ is degree $2$
  for $2 \leq i \leq \ell-1$ and either $u_1$ and $u_{\ell}$ are ot
  degree $2$ or $u_1=u_{\ell}$ ($O(n)$ steps). Set $\MARKED(u_i)=\TRUE$ for $2 \leq i \leq \ell-1$.

\item If $u_1=u_{\ell}$, $G$ is a cycle. Return ``YES''.

\item If $u_1 \neq u_{\ell}$ and there is an edge between $u_1$ and $u_{\ell}$, determine 
whether $G \setminus \{ u_2, \ldots, u_{\ell-1} \}$ is cyclically orientable. Return ``YES'' if and only if this subgraph is.

\item If $u_1 \neq u_{\ell}$ and no edge exists between $u_1$ and $u_{\ell}$, return to step (3).

\end{enumerate}

Since this algorithm takes $O(n)$ steps and then recurses to solving
the same problem on a smaller graph, its run time is $O(n^2)$. It is
easy to see that this algorithm is correct: after checking some base
cases, this algorithm looks for a cycle $C=(u_1, \ldots, u_{\ell},
u_1)$ such that $G=G' \cup_{e_1,e_2} C$. If no such cycle exists, the
algorithm returns ``NO'' -- which is correct according to theorem
\ref{Main}. Otherwise, assuming the correctness of the algorithm
inductively, the algorithm returns ``YES'' if and only if $G'$ is
cyclically orientable -- which is correct according to theorem \ref{Main}
again.

We now describe how to speed up this algorithm. The first key idea is that a cyclically orientable graph can not have too many edges:

\begin{prop}
Let $G$ be a cyclically orientable graph with $n$ vertices. Then $G$ has at most
 $2n-3$ edges.
\end{prop}

\begin{proof}
  We first prove this for $G$ two-connected. The bound is correct when
  $G$ is a single edge or a cycle. If $G=G' \cup_{e_1,e_2} Z$, where
  $G'$ is two-connected cyclically orientable graph on $n'$ vertices and $Z$ is 
a cycle with $z \geq 3$ vertices. If $G'$ has $e'$ edges, then $G$ has $e'+z-1$ 
edges and $n'+z-2$ vertices. We have (by induction on $n$) 
$$e'+z-1 \leq 2n'-3 +z -1 =2n'-z-4 \leq  2n'-2z-3=  2(n'+z-2)-3,$$
which is the desired result.

Now we do not assume that $G$ is two-connected. If $G_1$, \dots, $G_k$
are the two-connected components of $G$, with $G_i$ having $n_i$
vertioces and $e_i$ edges, then $G$ has at most $\sum n_i - (k-1)$ vertices
and has exactly $\sum e_i$ edges. We have
$$\sum e_i \leq \sum (2 n_i-3) =2 \sum n_i -3k \leq 2 \sum n_i - 2k+2 = 2 \left(
 \sum n_i - 
(k-1) \right)$$
which is the desired inequality.
\end{proof}

Therefore, if we begin ou algorithm by testing whether $G$ has more
than $2n-3$ edges or not, we can therefore assume that any computation
which runs in $O(e)$ time, where $e$ is the number of edges, in fact
runs in $O(n)$ time. We use two other ideas to speed up our
computation: First, suppose that $C$ is a chain in $G$. Then it is
clear that replacing $C$ by a chain with two edges does not effect the
cyclic orientability of $G$. Second, we remove the use of $\MARKED$ and
instead explore the vertices in an order which is automatically
non-redundant.

In detail, given a graph $G$, our algorithm first determines whether
$G$ has more than $2n-3$ edges. If so, it outputs no. If not, it then
breaks $G$ up into two-connected components. It then tests each two
connected component $G_i$, to see if it has more than $2 n_i-3$ edges,
where $n_i$ is the number of vertices of $G_i$. If any of them deos,
it outputs ``NO''. Once these preliminaries are done, it then carries
out the following computation on each component and returns ``YES'' if
and only if the following algorithm returns ``YES'' in every case.

Given a two-connected graph $G$ with $n$ vertices and $O(n)$ edges,
the following algorithm determines in $O(n)$ time whether this graph
is cyclically orientable. This algoirthm maintains a list $L$, which
should be thought of as a list of the degree two vertices that still
need to be dealt with. This list should be a doubly linked list of
pointers to vertices, with a pointer from each vertex back to its
point on the list and a pointer to the end of $L$ always maintained,
so that we can in time $O(1)$ delete elements from $L$, insert them
and find the end or beginning of $L$. Operations related to maintaining
the data structure $L$ will not be explicitly described.

\begin{enumerate}
  
\item Find all degree two vertices of $G$ and put them into a list $L$ as above.
  
\item If $L$ is empty, return ``NO''.

\item Take the first element $v$ of $L$. Find the unique path $u_1$,
  $u_2$, \dots, $u_{k-1}$, $u_k=v$, $u_{k+1}$, \dots, $u_{\ell}$ in
  $G$ such that $u_i$ is degree $2$ for $2 \leq i \leq \ell-1$ and
  either $u_1$ and $u_{\ell}$ are not degree $2$ or $u_1=u_{\ell}$. 

\item If $u_1=u_{\ell}$, $G$ is a cycle. Return ``YES''.
  
\item If there is an edge joining $u_1$ and $u_{\ell}$, delete the
  vertices $u_2$ through $u_{\ell-1}$ from $G$ and from $L$. If this
  make $u_1$ and/or $u_{\ell}$ have degree $2$, add them to the end of
  $L$. Return to step (2).
  
\item Delete the vertices $u_2$ through $u_{\ell-1}$ from $G$ and from
  $L$. Add a new vertex $w$ to $G$ with edges to $u_1$ and $u_{\ell}$.
  (Do NOT add $w$ to $L$.) Return to (2).

\end{enumerate}

A detailed analysis of correctness is left to the reader; this is
basically the same algorithm as before. To compute the running time,
note that step (1), which only occurs once, takes $O(n)$ steps. Call
an iteration of steps (3)-(6) ``type A'' if $\ell=3$ and there is no
edge joining $u_1$ to $u_3$ and ``type B'' otherwise.

The number of vertices of $G$ never increases. For a given type B run,
let $k$ be the decrease in the number of vertices of $G$; $k$ is
either $\ell-1$ or $\ell-2$ so $\ell \leq 3k$. Then that type B run
takes $O(\ell)=O(k)$ steps. In particular, the amount of time taken in
all type B runs is $O(n)$ as the sum of all the $k$'s must be at most
$n$. Also, there must be at most $n$ type B runs.

Only type B runs increase the size of $L$ and then only by one each
time. $L$ starts out with at most $n$ members and is increased by $1$
at most $n$ times during the run of the algorithm. Each type A run
decreases the size of $L$ by at least one. Thus, there are $O(n)$ type
A runs and each of these takes time $O(1)$. In total, at most $O(n)$
time is spent in type A runs, at most $O(n)$ time is spent in type B
runs and at most $O(n)$ time is spent in precomputation, so this algorithm runs in time $O(n)$.

\raggedright

\end{document}